\documentclass[11pt,twoside]{amsart}
\usepackage{graphicx}
\usepackage{color}
\usepackage{amssymb}
\usepackage{amsmath}
\usepackage{amsthm}
\usepackage{amsfonts}
\usepackage{amsmath, amsfonts, amssymb}
\newtheorem{theorem}{Theorem}[section]

\newtheorem{proposition}[theorem]{Proposition}

\newtheorem{remark}[theorem]{Remark}

\numberwithin{equation}{section}

\begin{document}
	\makeatletter
	\makeatother
	\small \begin{center}{\textit{ In the name of
				Allah, the Beneficent, the Merciful.}}\end{center}
	\large
			\hskip -0.2 cm
			\title{On equivalence of vector-valued maps}
				
	\author{Khadjiev Dj.$^{1}$, Bekbaev U.$^{2}$, Aripov R.$^{3}$}
	
	\thanks{{\scriptsize
			emails: $^1$khdjavvat@gmail.com; $^2$uralbekbaev@gmail.com; $^3$arrustamg@yandex.ru}}
	\maketitle
	\begin{center}
		{\scriptsize \address{$^1$National University of Uzbekistan, UzAS Institue of Mathematics, Tashkent, Uzbekistan}}\\
		{\scriptsize \address{$^2$Depart. of Mathematical and Natural Sciences, TTPU, Tashkent, Uzbekistan}}\\
		
	\end{center}

	\begin{abstract} An approach to the equivalence problem of vector valued maps is offered which, in particular, covers the equivalence problem of paths and patches of differential geometry with respect to different motion groups. In the last case, in contrary to differential geometry case, it does not need and does not use smoothness of paths and patches to get the corresponding main results.
		\end{abstract}
	
	\vskip 0.2 true cm

\section{Introduction}
The smooth paths and patches are research objects of differential geometry and their equivalence problem with respect to different motion groups are investigated by the use of differential calculus. In abstract differential field settings this kind problems are considered in \cite{Kh,B1,B2,B3} Of course, there are many other cases of "paths" and "patches" for which the equivalence problem can not be by handled by methods of differential geometry or even those methods may have no meaning in those other cases. In the current paper we are going to offer an approach (method) to solve the $G$-equivalence problem for system of vector-valued maps over any abstract underlying field, where $G$ is a subgroup of the corresponding general linear group. In general and special linear group cases we provide detailed results. One can see easily that the approach is natural, combines all similar problems and even in smooth path, patch cases it does not use their smoothness.

 In two dimensional real(rational) vector space case similar problems, with respect to orthogonal and special orthogonal motion groups, are considered in \cite{Besh,Kh1,Ash,Saf}.  
\section{ The main results}
Let $T$ be any fixed set, $n\geq 2$ be any natural number, $F$ be any field and $G\subset GL(n,F)$ a subgroup. We say that two maps $u: T\rightarrow F^n$, $v: T\rightarrow F^n$ are $G$-equivalent ($u(t)\simeq^Gv(t)$) if there exists $g\in G$ such that equality $v(t)=gu(t)$ holds true at all $t\in T$, where it is assumed that elements of $F^n$ are presented as column vectors. The set $\{gu(t): \ g\in G\}$ is said to be the $G$-orbit of $u(t)$. The rank of a map $u(t)$ we define as $\dim Span\{u(t): t\in T\}$ and denote by $rk(u)$. It is clear that $rk(u)\leq n$ and $rk(u)=k$ if and only if there exist $t_1, t_2,t_3,...,t_k\in T$ such that rank of matrix $[u(t_1), u(t_2),u(t_3),...,u(t_k)]$ is $k$ and $Span\{u(t): t\in T\}=Span\{u(t_i): i=1,2,3,..., k\}$. If $u(t)\simeq^Gv(t)$ then evidently $rk(u)=rk(v)$. If $A$ is a matrix we use notation $A_{i_1,i_2,...,i_k}$ for matrix consisting of the $i_1^{th}, i_2^{th},..., i_k^{th}$ rows of the matrix $A$, $I_k$ stands for the $k$-order identity matrix. Further it is assumed that $t_1, t_2,t_3,...,t_k\in T$ are fixed elements and we consider $G=GL(n,F)$-equivalence of maps $u(t)$ for which \[k=rk(u(t))=rk([u(t_1), u(t_2),u(t_3),...,u(t_k)])\], the set of such maps we denote by $V(t_1,t_2,...,t_k)$.

Proof of the following Proposition is evident.

\begin{proposition} If $rk([u(t_1), u(t_2),u(t_3),...,u(t_k)])=k$ then $rk(u(t))=k$ if and only if 
\[u(t)=\sum_{l=1}^k([u(t_1), u(t_2),u(t_3),...,u(t_k)]_{i_1,i_2,...,i_k})^{-1}u(t)_{i_1,i_2,...,i_k})_lu(t_l),\]  whenever
$\det[u(t_1), u(t_2),u(t_3),...\ \ ,u(t_k)]_{i_1,i_2,...,i_k}\neq 0$. In particular, if \\$\det[u(t_1), u(t_2),u(t_3),...,u(t_k)]_{j_1,j_2,...,j_k}\neq 0$ as well, then\[	
([u(t_1), u(t_2),u(t_3),...,u(t_k)]_{i_1,i_2,...,i_k})^{-1}u(t)_{i_1,i_2,...,i_k}=\] \[([u(t_1), u(t_2),u(t_3),...,u(t_k)]_{j_1,j_2,...,j_k})^{-1}u(t)_{j_1,j_2,...,j_k}\] at all $t\in T$.\end{proposition}

\begin{theorem}\label{1}  If $u(t), v(t)\in V(t_1,t_2,...,t_k)$ then $u(t)\simeq^Gv(t)$ if and only if there exist\\ $1\leq i_1<i_2<...<i_k\leq n,\ \ 1\leq j_1<j_2<...<j_k\leq n$ such that \[\det[u(t_1), u(t_2),u(t_3),...,u(t_k)]_{i_1,i_2,...,i_k}\neq 0,\ \ \det[v(t_1), v(t_2),v(t_3),...,v(t_k)]_{j_1,j_2,...,j_k}\neq 0\] and at all $t\in T$ the equality \[([u(t_1), u(t_2),u(t_3),...,u(t_k)]_{i_1,i_2,...,i_k})^{-1}u(t)_{i_1,i_2,...,i_k}=\]
	\[ ([v(t_1), v(t_2),v(t_3),...,v(t_k)]_{j_1,j_2,...,j_k})^{-1}v(t)_{j_1,j_2,...,j_k}\] holds true.\end{theorem}

\textbf{Proof.} Necessity. If $gu(t)=v(t)$ and $u(t)=\sum_{i=1}^k\alpha_k(t)u(t_i)$ then $v(t)=gu(t)=\sum_{i=1}^k\alpha_k(t)gu(t_i)=\sum_{i=1}^k\alpha_k(t)v(t_i)$ therefore there exist $1\leq i_1<i_2<...<i_k\leq n$, $1\leq j_1<j_2<...<j_k\leq n$ such that \[det[u(t_1), u(t_2),u(t_3),...,u(t_k)]_{i_1,i_2,...,i_k}\neq 0,\ \ det[v(t_1), v(t_2),v(t_3),...,v(t_k)]_{j_1,j_2,...,j_k}\neq 0\] and
\[([u(t_1), u(t_2),u(t_3),...,u(t_k)]_{i_1,i_2,...,i_k})^{-1}u(t)_{i_1,i_2,...,i_k}=\] 
 \[([v(t_1), v(t_2),v(t_3),...,v(t_k)]_{j_1,j_2,...,j_k})^{-1}v(t)_{j_1,j_2,...,j_k}\] 
 for all $t\in T$.

Sufficiency. There exists $g\in G$ such that \[ g[u(t_1), u(t_2),u(t_3),...,u(t_k)]= 
[v(t_1), v(t_2),v(t_3),...,v(t_k)].\] We show that in such case $gu(t)=v(t)$. Indeed
due to \[u(t)=(\sum_{l=1}^k([u(t_1), u(t_2),u(t_3),...,u(t_k)]_{i_1,i_2,...,i_k})^{-1}u(t)_{i_1,i_2,...,i_k})_lu(t_l)\]
one has \[gu(t)=\sum_{l=1}^k(([u(t_1), u(t_2),u(t_3),...,u(t_k)]_{i_1,i_2,...,i_k})^{-1}u(t)_{i_1,i_2,...,i_k})_lgu(t_l)=\]
\[\sum_{l=1}^k (([v(t_1), v(t_2),v(t_3),...,v(t_k)]_{j_1,j_2,...,j_k})^{-1}v(t)_{j_1,j_2,...,j_k})_lv(t_l)=v(t).\]

\begin{theorem} If $u^0_1(t), u^0_2(t),..., u^0_k(t)$ any functions such that \[[u^0(t_1),u^0(t_2),...,u^0(t_k)]=I_k\] then there exists unique, up to $G$-equivalence, map $u(t)\in V(t_1,t_2,...,t_k)$ such that 
	\[([u(t_1), u(t_2),u(t_3),...,u(t_k)]_{1,2,...,k})^{-1}u(t)_{1,2,...,k}=u^0(t)\] at all $t\in T$, where $u^0(t)$ stands for column vector with components $u^0_1(t), u^0_2(t),..., u^0_k(t)$.\end{theorem}

\textbf{Proof.} One can define the needed $u(t)\in V(t_1,t_2,...,t_k)$ in the following way: By definition $u_i(t)=u^0_i(t)$ whenever $1\leq i\leq k$ and $u_i(t)=0$ whenever $k<i\leq n$. It is easy to see that this map satisfies all requirements 
of the theorem. Uniqueness, up to $G$-equivalence, is an immediate consequence of Theorem \ref{1}. 

Let $F(x(t),x(t_1),x(t_2),...,x(t_k))$ stand for the field of rational functions over $F$ in column-vector variables $x(t)=(x_i(t))_{i=1,2,...,n}, x(t_j)=(x_i(t_j))_{i=1,2,...,n}$, where $j=1,2,...,k$. An element $f\in F(x(t),x(t_1),x(t_2),...,x(t_k))$ is said to be $G$-invariant if and only if \[f(gx(t),gx(t_1),gx(t_2),...,gx(t_k))=f(x(t),x(t_1),x(t_2),...,x(t_k))\] whenever $g\in G$. Further $F(x(t),x(t_1),x(t_2),...,x(t_k))^G $
stands for the set of all such $G$-invariants.

It is known that only in infinite field $F$ case any two rational functions over $F$ are formally equal if and only if their corresponding values over a set, where does not vanish some polynomial over $F$, are equal. As far as in the following theorem we need equality of rational functions as equality of their corresponding values we have to assume infinity of $F$.

\begin{theorem} If $F$ is infinite then the field $F(x(t),x(t_1),x(t_2),...,x(t_k))^G $ is generated by the system of components of the vector $([x(t_1), x(t_2),x(t_3),...,x(t_k)]^{-1})_{1,2,...,k}x(t)_{1,2,...,k}$ over $F$. Moreover the system is algebraic independent over $F$.\end{theorem} 

\textbf{Proof.} Let $f(x(t),x(t_1),x(t_2),...,x(t_k))\in F(x(t),x(t_1),x(t_2),...,x(t_k))^G $. Due to Theorem 2 the system of invariants $([x(t_1), x(t_2),x(t_3),...,x(t_k)]^{-1})_{1,2,...,k}x(t)_{1,2,...,k}$ defines  an unique, up to $G$-equivalence, map $x(t)$ therefore at least there exists a function $g(y_1,y_2,...,y_k)$ such that\\ $f(x(t),x(t_1),x(t_2),...,x(t_k))=g(([x(t_1), x(t_2),x(t_3),...,x(t_k)]^{-1})_{1,2,...,k}x(t)_{1,2,...,k})$, that is $f=g\circ h$, where \[h(x(t),x(t_1),x(t_2),...,x(t_k))= ([x(t_1), x(t_2),x(t_3),...,x(t_k)]^{-1})_{1,2,...,k}x(t)_{1,2,...,k}.\] To prove the rationality of $g$ it is enough to show that  the $G$-invariant rational map  $h$ has rational inverse. Note that the values of $h$ at different elements of a $G$-orbit are same. Therefore to define the value of the inverse function it is enough to give a representative of the corresponding orbit. For that let us consider the rational map $I: F^k\rightarrow F^{n+nk}$, where $I(y_1,y_2,...,y_k)=(y,e_1,e_2,...,e_k)$,  $y,e_1,e_2,...,e_k$ are vectors with $n$ components, all components of $y$ are zero except for the first $k$ components which are $y_1,y_2,...,y_k$, all components of $e_i$ are zero except for $i^{th}$ component which is $1$, $i=1,2,...,k$. It is easy to check that $I\circ h$ and $h\circ I$ are identical maps, so $g=f\circ I$ is a rational function. The algebraic independence of the system is evident. 

Instead of one map a family of maps $(u_s(t))_{s\in S}$ can be considered, where $S$ is a fixed set and one can define its rank as \[rk((u_s(t))_{s\in S})=\dim Span\{u_s(t): t\in T,\ s\in S\},\] fix some
$t_{1,s_1}, t_{2,s_2},..., t_{k,s_k}\in T$, where $s_i\in S$ at $i=1,2,...,k$, such that $Span\{u_{s_i}(t_{i,s_i}): i=1,2,...,k\}=$ $Span\{u_s(t): t\in T,\ s\in S\}$. Let $V(t_{1,s_1}, t_{2,s_2},..., t_{k,s_k})=$ \\ $\{(u_s(t))_{s\in S}: rk((u_s(t))_{s\in S})= rk[u_{s_1}(t_{1,s_1}), u_{s_2}(t_{2,s_2}),..., u_{s_k}(t_{k,s_k})]\}$.

Here are the formulations of the corresponding results in this case without proofs as far as proofs are similar to those of above.   

\begin{theorem} If $(u_s(t))_{s\in S}, (v_s(t))_{s\in S}\in V(t_{1,s_1}, t_{2,s_2},..., t_{k,s_k})$ then\\ $(u_s(t))_{s\in S}\simeq^G(v_s(t))_{s\in S}$, that is $gu_s(t)=v_s(t)$ for some $g\in G$ at all $s\in S, t\in T$, if and only if there exist $1\leq i_1< i_2< ... < i_k\leq n$, $1\leq j_1< j_2< ... < j_k\leq n$
	such that \[\det([u_{s_1}(t_{1,s_1}), u_{s_2}(t_{2,s_2}),..., u_{s_k}(t_{k,s_k})]_{i_1,i_2,...,i_k})\neq 0,\] \[\det([v_{s_1}(t_{1,s_1}), v_{s_2}(t_{2,s_2}),..., v_{s_k}(t_{k,s_k})]_{j_1,j_2,...,j_k})\neq 0\] and \[([u_{s_1}(t_{1,s_1}), u_{s_2}(t_{2,s_2}),..., u_{s_k}(t_{k,s_k})]_{i_1,i_2,...,i_k}^{-1}u_s(t)_{i_1,i_2,...,i_k}=\]
	\[([v_{s_1}(t_{1,s_1}), v_{s_2}(t_{2,s_2}),..., v_{s_k}(t_{k,s_k})]_{j_1,j_2,...,j_k}^{-1}v_s(t)_{j_1,j_2,...,j_k}\] for all $t\in T$ and $s\in S$.\end{theorem}

\begin{theorem} If $(u^0_{i,s}(t))_{s\in S, i=1,2,...,k}$ is any system of functions such that \[[u^0_{s_1}(t_{1,s_1}),u^0_{s_2}(t_{2,s_2}),...,u^0_{s_k}(t_{k,s_k})]=I_k\] then there exists unique, up to $G$-equivalence, system of maps $(u_s(t))_{s\in S}\in V(t_{1,s_1}, t_{2,s_2},..., t_{k,s_k})$ such that 
	\[([u_{s_1}(t_{1,s_1}), u_{s_2}(t_{2,s_2}),..., u_{s_n}(t_{k,s_k})]_{1,2,...,k}^{-1}u_s(t)_{1,2,...,k}=
	u^0_s(t)\] at all $t\in T$, where $u^0_s(t)$ stands for column vector with components $u^0_{1,s}(t), u^0_{2,s}(t),..., u^0_{k,s}(t)$.\end{theorem}

Let $F((x_s(t))_{s\in S},x_{s_1}(t_{1,s_1}),x_{s_2}(t_{2,s_2}),...,x_{s_k}(t_{k,s_k}))$ stand for the field of rational functions over $F$ in column-vector variables $x_s(t)=(x_{i,s}(t))_{i=1,2,...,n}, x_{s_j}(t_{j,s_j})=(x_{i}(t_{j,s_j}))_{i=1,2,...,n}$, where $j=1,2,...,k$. An element  $f\in  F((x_s(t))_{s\in S},x_{s_1}(t_{1,s_1}),x_{s_2}(t_{2,s_2}),...,x_{s_k}(t_{k,s_k}))$ is said to be $G$-invariant if and only if \[f((gx_s(t))_{s\in S},gx_{s_1}(t_{1,s_1}),gx_{s_2}(t_{2,s_2}),...,gx_{s_k}(t_{k,s_k}))=f((x_s(t))_{s\in S},x_{s_1}(t_{1,s_1}),x_{s_2}(t_{2,s_2}),...,x_{s_k}(t_{k,s_k}))\] whenever $g\in G$. Further  $F((x_s(t))_{s\in S},x_{s_1}(t_{1,s_1}),x_{s_2}(t_{2,s_2}),...,x_{s_k}(t_{k,s_k}))^G $
stands for the set of all such $G$-invariants.
\begin{theorem} If $F$ is infinite then the field $F((x_s(t))_{s\in S},x_{s_1}(t_{1,s_1}),x_{s_2}(t_{2,s_2}),...,x_{s_k}(t_{k,s_k}))^G$ is generated by the system of components of the vectors \[([x_{s_1}(t_{1,s_1}), x_{s_2}(t_{2,s_2}),..., x_{s_k}(t_{k,s_k})]_{1,2,...,k})^{-1}x_s(t)_{1,2,...,k},\ \mbox{where}\ s\in S\] over $F$. Moreover the system is algebraic independent over $F$.\end{theorem} 

What about similar results when one considers instead of $GL(n,F)$ some subgroup $G$ of it? In $k=n$ case the following results are valid.  

\begin{theorem} If $(u_s(t))_{s\in S}, (v_s(t))_{s\in S}\in V(t_{1,s_1}, t_{2,s_2},..., t_{n,s_n})$ then\\ $(u_s(t))_{s\in S}\simeq^G(v_s(t))_{s\in S}$ if and only if
  \[([u_{s_1}(t_{1,s_1}), u_{s_2}(t_{2,s_2}),..., u_{s_n}(t_{n,s_n})]^{-1}u_s(t)=([v_{s_1}(t_{1,s_1}), v_{s_2}(t_{2,s_2}),..., v_{s_n}(t_{n,s_n})]^{-1}v_s(t),\] whenever $t\in T$, $s\in S$, and \[[v_{s_1}(t_{1,s_1}), v_{s_2}(t_{2,s_2}),..., v_{s_n}(t_{n,s_n})]([u_{s_1}(t_{1,s_1}), u_{s_2}(t_{2,s_2}),..., u_{s_n}(t_{n,s_n})]^{-1}\in G.\]\end{theorem}

Note that for some $G$, for example for $G=SL(n,F)$, the condition $g_2g_1^{-1}\in G$ is equivalent to a system of equalities $(f_p(g_1)=f_p(g_2))_{p\in P}$ ($\det(g_1)=\det(g_2)$ if $G=SL(n,F) $) and in such case the last condition of the above theorem can be replaced by 
\[(f_p([u_{s_1}(t_{1,s_1}), u_{s_2}(t_{2,s_2}),..., u_{s_n}(t_{n,s_n})])=f_p([v_{s_1}(t_{1,s_1}), v_{s_2}(t_{2,s_2}),..., v_{s_n}(t_{n,s_n})]))_{p\in P}\] 
\[(\det[u_{s_1}(t_{1,s_1}), u_{s_2}(t_{2,s_2}),..., u_{s_n}(t_{n,s_n})]=\det[v_{s_1}(t_{1,s_1}), v_{s_2}(t_{2,s_2}),..., v_{s_n}(t_{n,s_n})]\ \mbox{if}\ G=SL(n,F)).\]

\begin{theorem} If $(u^0_{i,s}(t))_{s\in S, i=1,2,...,n}$ is any system of functions such that \[[u^0_{s_1}(t_{1,s_1}),u^0_{s_2}(t_{2,s_2}),...,u^0_{s_k}(t_{n,s_n})]=I_n\] there exists unique, up to $G$-equivalence, system of maps $(u_s(t))_{s\in S}\in V(t_{1,s_1}, t_{2,s_2},..., t_{n,s_n})$ such that $[u_{s_1}(t_{1,s_1}), u_{s_2}(t_{2,s_2}),..., u_{s_n}(t_{n,s_n})]\in G$ and
	\[([u_{s_1}(t_{1,s_1}), u_{s_2}(t_{2,s_2}),..., u_{s_n}(t_{n,s_n})]^{-1}u_s(t)=
	u^0_s(t)\] at all $t\in T$, where $u^0_s(t)$ stands for column vector with components $u^0_{1,s}(t), u^0_{2,s}(t),..., u^0_{n,s}(t)$.\end{theorem}

In this theorem the condition $[u_{s_1}(t_{1,s_1}), u_{s_2}(t_{2,s_2}),..., u_{s_n}(t_{n,s_n})]\in G$ merely means\\  $\det[u_{s_1}(t_{1,s_1}), u_{s_2}(t_{2,s_2}),..., u_{s_n}(t_{n,s_n})]=1$ if $G=SL(n,F)$. 

\begin{theorem} If $F$ is infinite then  \[F((x_s(t))_{s\in S},x_{s_1}(t_{1,s_1}),x_{s_2}(t_{2,s_2}),...,x_{s_n}(t_{n,s_n}))^G=\overline{F}(x_{s_1}(t_{1,s_1}),x_{s_2}(t_{2,s_2}),...,x_{s_n}(t_{n,s_n}))^G,\] and the system of entries of vector variables $x_{s_1}(t_{1,s_1}),x_{s_2}(t_{2,s_2}),...,x_{s_n}(t_{n,s_n})$ is algebraic independent over the field $\overline{F}$, where $\overline{F}=F(([x_{s_1}(t_{1,s_1}), x_{s_2}(t_{2,s_2}),..., x_{s_n}(t_{n,s_n})]^{-1}x_s(t))_{s\in S}).$\end{theorem} 

An evident advantage of the last theorem is that it reduces description of the field of $G$-invariant rational functions over $F$ of system of maps to description of the field of $G$-invariant rational functions over another field of finite system, cardinality $n$, of vectors. It should be noted also that $G$-invariants of finite system of vectors is relatively well studied part of the invariant theory \cite{H}.  

 In $G=SL(n,F)$ case the last theorem says that \[F((x_s(t))_{s\in S},x_{s_1}(t_{1,s_1}),x_{s_2}(t_{2,s_2}),...,x_{s_n}(t_{n,s_n}))^G\] is generated over $F$ by the system of entries of $([x_{s_1}(t_{1,s_1}), x_{s_2}(t_{2,s_2}),..., x_{s_n}(t_{n,s_n})]^{-1}x_s(t))_{s\in S}$ and 
 $\det [x_{s_1}(t_{1,s_1}), x_{s_2}(t_{2,s_2}),..., x_{s_n}(t_{n,s_n})]$.
 
 Note that if $G=SL(n,F)$ and $k< n$ then $G$- equivalence of system of vector-valued maps is reduced to the case of $G=GL(n,F)$.

 \begin{remark} If $G\subset GL(n,F)$ is a fixed subgroup it is clear that $gu_s(t)+u_0=v_s(t)$ for all $t\in T, s\in S$ and some $g\in G$, $u_0\in F^n$ if and only if $g(u_s(t)-u_{s_0}(t_0))=v(t)-v_{s_0}(t_0)$ and $gu_{s_0}(t_0)+u_0=v_{s_0}(t_0)$, where $s_0\in S$, $t_0\in T$ are fixed elements. Therefore $G_{aff}=F^n\rtimes G$- equivalence problem of system of vector-valued maps can be solved as far as $u(t)\simeq^{G_{aff}}v(t)$ if and only if $(u_s(t)-u_{s_0}(t_0))_{s\in S}\simeq^{G}(v_s(t)-v_{s_0}(t_0))_{s\in S}.$\end{remark}

\end{document}